# Optimizing the Geometry of an L-Shaped Building to Enhance Energy Efficiency and Sustainability

Ewa ROKITA-MAGDZIARZ[1], Barbara GRONOSTAJSKA[2], Marcin MAGDZIARZ[3]*

[1]*Rokita-Project Architectural Office, Forsycjowa 7, Wrocław, 51-253, Poland*

[2]*Department of Architectural and Construction Design, Faculty of Architecture, Wrocław University of Science and Technology, Wyspianskiego 27, Wrocław, 50-370, Poland*

[3]*Hugo Steinhaus Center, Department of Applied Mathematics, Faculty of Pure and Applied Mathematics, Wrocław University of Science and Technology, Wyspianskiego 27, Wrocław, Poland*

*Abstract –* **The geometric form of a building strongly influences its material use, heat losses, and energy efficiency. This paper presents an analytical optimization of L-shaped residential buildings aimed at minimizing the external surface area for a prescribed volume. Both symmetric and asymmetric configurations are examined under realistic design constraints, including fixed or bounded wing aspect ratios and fixed building height. Using explicit optimization methods and Karush-Kuhn-Tucker conditions, closed-form expressions for the optimal geometric parameters and minimal envelope area are derived. The results show that unconstrained optimization leads to degenerate cuboid shapes, highlighting the importance of geometric constraints to preserve the L-shaped form. The obtained results provide practical design guidelines for architects and engineers, supporting informed early stage decisions that balance functional requirements, regulatory constraints, architectural intent, and energy performance. Case studies of existing houses demonstrate that the proposed approach can reduce external surface area or confirm near-optimality of practical designs, supporting energy-efficient early-stage architectural decisions.**

---

* Corresponding author.
E-mail address: marcin.magdziarz@pwr.edu.pl

## 1. INTRODUCTION

The building geometry is one of the most influential determinants of its structural behavior, material consumption, and energy requirements. Within the broad spectrum of residential typologies, the L-shaped house stands out as a particularly significant example in modern architecture. L-shaped houses represent a commonly used residential layout and have gained increasing interest from investors owing to their functional flexibility and architectural attractiveness. According to the National Association of Home Builders annual survey, more than 18% of newly requested residential floor plans in the early weeks of 2024 featured an L-shaped configuration [1]. Composed of two orthogonally intersecting wings, this form provides both practical and aesthetic benefits in traditional as well as contemporary designs. By disrupting continuous wind flow around the building, the L-shaped geometry enhances wind resistance and structural stability when compared with more linear building forms. Despite this, their design is frequently driven by tradition, stylistic norms, or regulatory constraints rather than by analytical optimization. Consequently, many realized examples diverge considerably from theoretically optimal solutions, resulting in unnecessary increases in material demand, construction costs, and operational energy use.

The relationship between building form and energy performance has been extensively documented in the literature. Depecker et al. [2] demonstrated the impact of building shape on heating demand, showing that compact envelopes reduce heat losses. This principle was further emphasized by Hegger et al. [3] in their handbook on sustainable architectural design, underscoring the importance of geometry in lowering operational energy consumption. In a comprehensive life-cycle energy assessment, Ramesh et al. [4] showed that parameters related to form can represent a significant share of total building energy use. Numerous other studies have sought to model and quantify this relationship, including those by Catalina et al. [5], Ourghi et al. [6], Schlueter and Thesseling [7], and Jedrzejuk [8, 9], all of which confirmed the strong link between compactness and performance. Steemers [10] further argued that decisions made at the morphological level of architectural design are among the most effective levers for reducing energy demand.

More recently, D'Amico and Pomponi [11] proposed a dimensionless compactness index that enables scale-independent assessment of building geometries, addressing the shortcomings of conventional surface-to-volume ratios. Their findings underline that geometry is as critical to sustainability as material selection or building systems.

Studies at the urban and typological scales corroborate these conclusions. Steadman's investigations into design evolution and building typologies [12, 13] illustrate the persistent influence of form on efficiency, while Julia et al. [14] and Hargreaves [15] incorporated typological and density-related variables into urban energy models. Angel's global analysis of urban expansion [16], together with United Nations reports [17] and policy instruments such as the London Strategy [18] and California's energy regulations [19], emphasize that residential design choices have enduring consequences for sustainable urban development.

An additional foundation for this research is provided by the theoretical work on architectural form generation. Early studies [20-22] introduced shape grammars, demonstrating that design processes can be formalized through rule-based systems. Mitchell [23] strengthened this perspective by integrating formal logic into computational design, while Knight [24] developed these ideas into a comprehensive framework for systematic shape exploration. Synthesized by Oxman [25], these contributions form the basis of much contemporary computational design theory, which conceives architecture as the outcome of generative rules informed by performance criteria. Within this paradigm, geometry becomes an explicit and controllable design variable, subject to optimization alongside other constraints.

Optimization has consequently emerged as a central component of sustainable architectural design. Convex optimization, in particular, provides a rigorous mathematical framework for addressing constrained design problems and identifying global optima in form-finding applications [26]. Okeil [27], Caruso et al. [28], Jin and Jeong [29], and Vartholomaios [30] employed such techniques in architectural contexts, spanning parametric design studies and performance-driven optimization. Hachem [31, 32] extended these methods to the neighborhood scale, optimizing residential morphologies to enhance solar access and energy efficiency. Collectively, these works demonstrate that optimization techniques constitute practical and effective tools for advancing sustainability in the built environment.

Building on this body of knowledge, this paper introduces a rigorous mathematical framework for optimizing L-shaped geometry houses. In what follows, we develop and apply a rigorous optimization framework for the geometry of L -shaped houses. Section 2 presents an analytical formulation of the building envelope area and derives optimal geometric configurations for both symmetric and asymmetric L-shaped layouts under a range of realistic design scenarios, including fixed volume, fixed or bounded wing aspect ratios, and prescribed building height. Closed form solutions are obtained using classical optimization techniques and Karush-Kuhn-Tucker conditions, providing clear design rules for minimizing the external surface area. Section 3 illustrates the practical applicability of the theoretical results through case studies of existing residential designs, demonstrating how the proposed framework can inform early-stage architectural decisions aimed at improving energy efficiency and sustainability. Section 4 concludes the paper by summarizing the key results, discussing their implications for sustainable housing design, and suggesting directions for future research.

## 2. OPTIMIZATION OF GEOMETRIC DIMENSIONS

L-shaped houses are a widely adopted residential layout and have become increasingly attractive to investors due to their functional versatility and architectural appeal. In the early weeks of 2024, more than 18% of newly requested residential floor plans in the National Association of Home Builders annual survey incorporated an L-shaped configuration [1]. Characterized by two intersecting wings, this configuration offers both practical and aesthetic advantages in traditional and contemporary designs. The L-shaped form improves wind resistance by breaking up wind flow across the structure, contributing to enhanced stability compared to more linear building layouts. It also facilitates effective drainage of rain, snow, and ice by allowing roof planes to be oriented in multiple directions, thereby reducing the likelihood of water accumulation and structural deterioration. The compact yet flexible geometry supports efficient zoning of interior spaces, promoting natural ventilation and daylight penetration, which can contribute to reduced energy demands for heating and cooling when combined with appropriate materials and construction methods. From a design perspective, L-shaped houses are valued for their balanced proportions and adaptability, allowing seamless integration of outdoor courtyards, patios, or extensions. Additionally, the configuration can accommodate attic or loft spaces and enables the incorporation of dormers or future expansions. Due to their resilience, energy efficiency potential, and design flexibility, L-shaped houses are often associated with increased property value and sustained investor interest.

Considering the growing popularity of L-shaped houses, determining their optimal forms is an important and timely task.

## 2.1. Symmetric L-Shaped House

An illustrative sketch of a symmetric L-shaped house is shown in Fig. 1. The building footprint is composed of two perpendicular rectangular wings that intersect orthogonally to form a corner typical of L-plans. Both wings share the same geometric parameters: plan length L (measured along the longitudinal axis of each wing), uniform wing width B (measured across the short axis), and constant height H (vertical dimension).

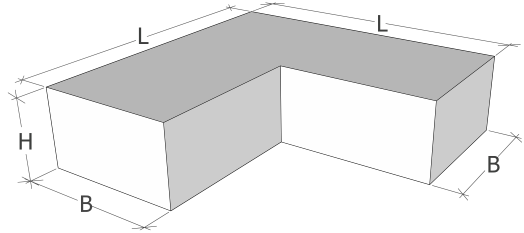

Fig. 1. View of a typical symmetric L-shaped house.

The volume of the house is given by:
$$V = H(2LB - B^2). \tag{1}$$

The external surface area of the L-shaped house, also called the building envelope), is equal to:
$$S = 4LH + 2LB - B^2. \tag{2}$$

Note that the first component represent the area of the walls. The remaining part is equal to the area of the roof. It should be noted that the area of the building in contact with the ground is not included in $S$.

The main aim of this section is to identify the optimal L-shaped house geometry that minimizes the external surface area $S$. We will examine different scenarios.

### 2.1.1 First Case: Fixed Volume V

In this scenario, the building volume V is assumed to be fixed (for example as specified by the investor). The objective is to identify the values of L, B and H that minimize the building's external surface area S.

However, some standard calculations show that in such setting the minimal area of S is obtained for L=B. Then, the house takes on a cuboid shape, which means that the L-form of the building is lost. This is therefore not a satisfactory solution.

A similar situation occurs when we assume that the floor area F=2LB-B^2 and height H of the building are fixed. Then the minimal envelope area is also obtained for L=B. Again, this solution yields a cuboid-shaped configuration, which is not desirable in this context.

*2.1.2 Second Case: Fixed Wing Aspect Ratio $r = L/B$*

At the preliminary stage of the design process, it is justifiable for an investor or designer to treat the wing aspect ratio r=L/B as a prescribed parameter, since this ratio encapsulates functional, regulatory, and contextual constraints that typically precede detailed geometric optimization. The wing aspect ratio determines the proportional relationships of interior spaces, thereby affecting room layouts, circulation efficiency, and furniture arrangement, and is consequently closely linked to programmatic requirements defined by the investor. By fixing r, the designer can ensure compliance with established spatial standards, including minimum room widths, corridor dimensions, and modular construction systems.

From a regulatory and site-planning standpoint, zoning regulations, setback requirements, and plot geometry may impose preferred or restrictive proportions on building wings. Prescribing the aspect ratio at an early design stage facilitates regulatory compliance while reducing the complexity of subsequent optimization procedures. Furthermore, construction-related factors-such as structural span limitations, standardization of prefabricated components, and cost-efficiency considerations - often favor specific length-to-width ratios, rendering a fixed r both practical and economically justified.

Finally, considerations of architectural coherence and design intent may also warrant fixing the wing aspect ratio in order to preserve visual balance, façade rhythm, and compatibility with the surrounding built environment.

Assuming that V and r=L/B are fixed, we get from (1)-(2) that

$$H = \frac{V}{2LB - B^2}$$

and

$$S = S(B) = \frac{4Vr}{B(2r - 1)} + B^2(2r - 1) \tag{3}$$

Therefore, the derivative of $S$ with respect to $B$ equals

$$\frac{dS}{dB} = 2B(2r-1) - \frac{4Vr}{B^2(2r-1)}$$

Solving equation $\frac{dS}{dB} = 0$ we get that the minimum of $S(B)$ is obtained for

$$B_{min} = \left(\frac{2Vr}{4r^2 - 4r + 1}\right)^{1/3} \tag{4}$$

Consequently, the remaining optimal parameters are

$$L_{min} = rB_{min} = r\left(\frac{2Vr}{4r^2 - 4r + 1}\right)^{1/3} \tag{5}$$

$$H_{min} = \frac{V}{2L_{min}B_{min} - B_{min}^2} = \frac{V}{(2r-1)\left(\frac{2Vr}{(2r-1)^2}\right)^{2/3}} = \frac{[(2r-1)V]^{1/3}}{(2r)^{2/3}}. \tag{6}$$

The minimal envelope is equal to

$$S_{min} = S(B_{min}) = 3\left(\frac{4V^2r^2}{2r-1}\right)^{1/3} \tag{7}$$

see Eq. (3).

The expressions given in (4)-(7) offer a practical basis for the design of a L-shaped house that attains the minimum external surface area $S_{min}$ for a prescribed volume $V$ and ratio $r$. Utilizing these relations enables designers to optimize the building geometry, thereby reducing material consumption and construction costs. In addition, minimizing the external surface area enhances energy efficiency by decreasing heat losses through the building envelope, which in turn improves overall sustainability and lowers long-term operational costs.

Let us analyze the following particular example: take $V = 300 m^3$ and $r = 2$. Then, from (4)-(7) we obrtain the following optimal parameters: $B_{min} = 5{,}11m, L_{min} = 10{,}22m, H_{min} = 3{,}83m, S_{min} = 234{,}89m^2$.

It is worth noting that, from a practical perspective, the parameters derived above are well suited to real-world residential construction. They provide viable design options while ensuring structural efficiency, economic feasibility, and improved energy performance. Consequently, these results offer a dependable basis for architects and engineers aiming to implement optimized house designs in practice.

Fig. 2 illustrates the external surface area $S$ of a symmetric L-shaped house as a function of the wing width $B$ and the wing aspect ratio $r = L/B$ for a fixed building volume $V = 300 \text{ m}^3$. The surface is presented as a three-dimensional plot, where each curve corresponding to a fixed value of $r$ exhibits a clear minimum. These minima, marked by red dots, identify the optimal wing width $B_{min}$ that minimizes the external surface area $S$ for a given aspect ratio. The figure visually demonstrates that, although $S$ increases for very small or very large values of $B$, an intermediate optimal configuration exists for each $r$, thereby confirming the analytical results derived from Eq. (3) and highlighting the dependence of geometric optimization on the prescribed aspect ratio.

Compactness measure. A common method for characterizing the relationship between a building's external surface area $S$ and its volume $V$ is the surface-to-volume ratio [2,3]. Lower values of this ratio are typically associated with greater geometric compactness. However, an important drawback of the $S/V$ ratio is its dependence on scale: geometrically similar buildings of different sizes (for example, cubes with different edge lengths) yield different $S/V$ values.

To overcome this limitation, a recently introduced compactness metric has been proposed, which evaluates building compactness purely on the basis of shape, independent of absolute size [11]. This measure is defined as $\frac{S}{S_{min}}$, see [11]. Here $S$ denotes the external surface area of the building and $S_{min}$ is the minimum external surface area required to enclose a given volume $V$. This metric is particularly useful in the early design phase, as it quantitatively indicates how closely a given form approaches optimal compactness, thereby supporting informed design decisions and the exploration of more efficient geometric alternatives.

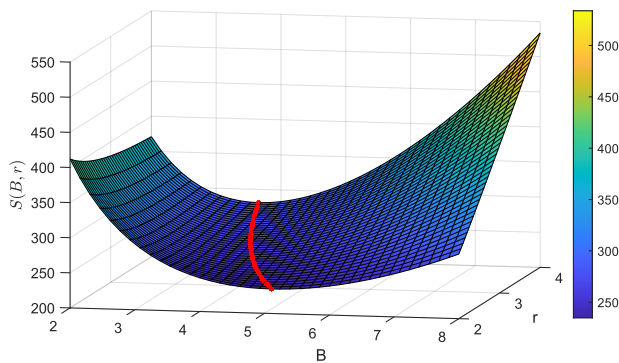

Fig. 2. External surface $S = S(B, r)$ as a function of parameters $S$ and $r$, see Eq. (3). Here $V = 300 m^3$. For each fixed $r$, there is a minimium of $S$ attained at $B_{min}$ and denoted here by red dot.

We now evaluate this compactness measure for the L-shaped house with fixed $V$ and $r$. Substituting equations (3) and (7) yields

$$\frac{S}{S_{\min}} = \frac{\frac{4Vr}{B(2r-1)} + B^2(2r-1)}{3\left(\frac{4V^2r^2}{2r-1}\right)^{1/3}} = \frac{4Vr + B^3(2r-1)^2}{3 \cdot B[2Vr(2r-1)]^{2/3}} \quad (8)$$

By construction, $S/S_{\min} \geq 1$, with equality attained for the optimal values (4)-(6). Consequently, this compactness measure provides a quantitative indicator of the geometric efficiency of the house design. It serves as a valuable tool for assessing shape optimality, guiding design refinement, and identifying opportunities for reducing both energy demand and construction costs.

*2.1.3 Third Case: Wing Aspect Ratio r from a Fixed Range*

In the early design stage it is reasonable for an investor or architect to assume that the wing aspect ratio $r = L/B$ belongs to a prescribed range, $r \in [a, b]$, rather than being completely free. The main reasons are practical and come from constraints that typically precede any geometric optimization of the envelope. The ratio $r$ strongly influences the usability of each wing: room proportions, corridor efficiency, furniture layouts, and adjacency requirements (e.g., kitchen-dining-living sequences, bedroom zoning) all impose limits on how narrow or elongated a wing can be. Very large $r$ tends to produce "corridor-like" plans with inefficient circulation, while very small $r$ pushes the wing toward a near-square block and can undermine the intended L-configuration. Moreover, plot geometry, setbacks, easements, required distances to boundaries, and local planning rules often restrict feasible building proportions. For example, a narrow site may cap $B$, while boundary offsets and courtyard requirements may cap $L$. Treating $r$ as lying in a feasible interval ensures compliance from the outset and avoids exploring geometries that cannot be permitted or built on the given parcel. Additionally, spans, roof framing, and load-bearing layouts tend to favor certain widths $B$ (linked to economical joist/truss spans and standard structural modules). Meanwhile, excessive length $L$ can increase expansion joints, lateral stability demands, and foundation complexity. A bounded range for $r$ reflects what is structurally straightforward and cost-effective to construct with typical systems.

Although optimization may later minimize surface area for a given volume, extreme aspect ratios can increase envelope area, thermal bridging risk, and heat loss, and may complicate airtightness and insulation continuity at corners and junctions. Constraining $r$ to a realistic interval helps keep early concepts aligned with energy-efficiency targets. The L-shaped topology is often chosen for a clear compositional reason (courtyard formation, façade rhythm, massing balance). If $r$ drifts too far, the design can either collapse toward a compact block (losing the

L-shape) or become overly stretched and visually unbalanced. A predefined range preserves the intended character while still allowing optimization within meaningful design freedom. For these reasons, specifying an admissible interval $r \in [a, b]$ early is a practical way to encode programmatic, regulatory, structural, energy, and aesthetic constraints-reducing the design space to solutions that are both feasible and architecturally relevant.

In this scenario, we assume that the ratio $r$ is constrained to a prescribed interval, $r \in [a, b]$. As before, the volume $V$ is given. Under these assumptions, our objective is to determine the optimal house geometry that minimizes the external surface area $S$, as defined in (3).

To this end, we employ the method of Lagrange multipliers combined with the Karush-Kuhn-Tucker (KKT) conditions [26]. Specifically, we seek the optimal parameters $(B^*, r^*)$ that minimize the external surface area, in accordance with (3).

More precisely, we seek to find the minimum of the function

$$S(B, r) = \frac{4Vr}{B(2r - 1)} + B^2(2r - 1)$$

with constraints $a \leq r \leq b$ and $B > 0$. The constants $V, a$, and $b$ are fixed, with $V > 0, a > 1, b > 1$ and $a < b$. The constraint functions are:

$$g_1(r) = a - r \leq 0 \text{ (i.e., } r \geq a),$$
$$g_2(r) = r - b \leq 0 \text{ (i.e., } r \leq b).$$

The corresponding Lagrangian $\mathcal{L} = \mathcal{L}(B, r, \lambda_1, \lambda_2)$ is defined as:

$$\mathcal{L} = \frac{4Vr}{B(2r-1)} + B^2(2r - 1) + \lambda_1(a - r) + \lambda_2(r - b).$$

Here $\lambda_1$ and $\lambda_2$ are the Lagrange multipliers.

The necessary conditions for a local minimum at $(B^*, r^*)$ are the KKT conditions:

I. Stationarity conditions $\frac{\partial \mathcal{L}}{\partial B} = 0$ and $\frac{\partial \mathcal{L}}{\partial r} = 0$.

The first condition yields

$$B = \sqrt[3]{\frac{2Vr}{(2r-1)^2}},$$

the second one gives

$$-\frac{4V}{B(2r-1)^2} + 2B^2 - \lambda_1 + \lambda_2 = 0.$$

II. Primal Feasibility: The constraints must be satisfied

$$a - r \leq 0, r - b \leq 0.$$

III. Dual Feasibility: The Lagrange multipliers must be non-negative

$$\lambda_1 \geq 0, \lambda_2 \geq 0.$$

IV. Complementary Slackness:

$$\lambda_1(a - r) = 0, \lambda_2(r - b) = 0.$$

Solving the above KKT system, we obtain that the unique points $(B^*, r^*)$ corresponding to the minimum of $S$ with given constrains, are equal to

$$r^* = a \tag{9}$$

and

$$B^* = \left(\frac{2Va}{(2a-1)^2}\right)^{1/3}. \tag{10}$$

The corresponding minimal value $S^* = S(B^*, r^*)$ is therefore equal to

$$S^* = 3\left(\frac{4V^2 a^2}{2a-1}\right)^{\frac{1}{3}}. \tag{11}$$

The remaining optimal parameters are

$$L^* = aB^* = a\left(\frac{2Va}{4a^2 - 4a + 1}\right)^{1/3}, \tag{12}$$

$$H^* = \frac{V}{2L^* B^* - (B^*)^2} = \frac{[(2a-1)V]^{1/3}}{(2a)^{2/3}}. \tag{13}$$

The above results (9)-(13) give a practical basis for the design of a L-shaped house that attains the minimum external surface area for a prescribed volume and ratio $r$ from a predefined interval $[a, b]$. They also show that the optimal value of the ratio $r \in [a, b]$ is equal to $a$, i.e. the left end of the interval $[a, b]$.

Let us consider the following example. Take $V = 200 m^3, a = 3$ and $b = 4$. Then, from (9)-(13) we get that the optimal parameters are: $r^* = 3, B^* = 3{,}63 m, L^* = 10{,}90 m, H^* = 3{,}03 m$ and $S^* = 198{,}12 m^2$. Figure 3 presents the corresponding external surface area

$S(B,r)$. The surface is shown as a three-dimensional plot over the admissible range of $r$ constrained by $a = 3$ and $b = 4$. The visualization clearly illustrates how the external surface varies within this domain and highlights the influence of both geometric parameters on the building envelope. The global minimum of $S$ corresponding to the optimal solution obtained analytically from the constrained optimization problem, is marked by a red dot. This minimum occurs at the lower bound of the admissible aspect-ratio interval, confirming that the optimal configuration is attained for $r^* = a$. The figure thus provides a graphical validation of the KKT-based analysis and illustrates the effect of imposing bounds on the wing aspect ratio during the design process.

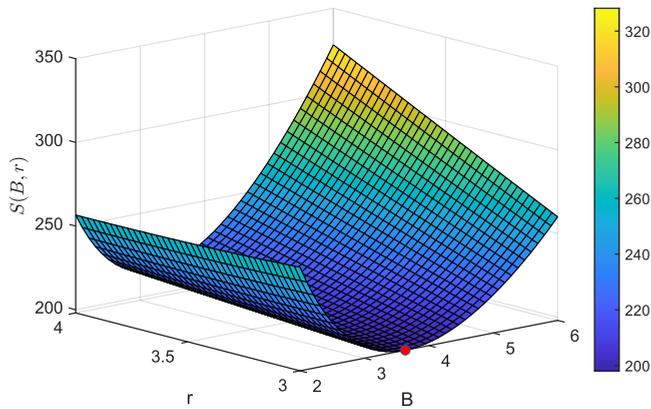

Fig. 3. Plots of the external surface $S(B,r)$ with $V = 200m^3$. For constraints $a = 3, b = 4$ the global minimal surface derived using formula (11) is equal to $198,12m^2$ and is denoted by the red dot.

## 2.2. Asymmetric L-Shaped House

In this section we will find optimality condition for an asymmetric L-shaped house. Illustrative sketch of the asymmetric L-shaped house is shown in Fig. 4. It is composed of two perpendicular wings of unequal dimensions. Unlike the symmetric configuration, the wings differ in length and/or width, resulting in a non-uniform footprint that departs from mirror symmetry. The intersection of the two rectangular wings forms a characteristic L-shaped plan, creating a sheltered inner corner that may be used as a courtyard, terrace, or protected outdoor space.

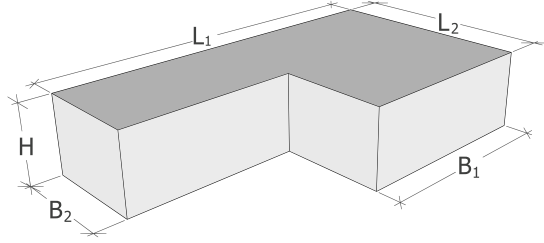

Fig. 4. View of an asymmetric L-shaped house.

The variation in wing proportions allows the building to respond more effectively to site-specific constraints such as plot geometry, orientation, access points, or solar exposure. The asymmetric L-shaped configuration offers greater design freedom than its symmetric counterpart, making it particularly suitable for irregular plots and context-driven residential designs, while retaining the climatic and spatial advantages associated with L-shaped forms.

Let us introduce the wing aspect ratios

$$r_1 = \frac{B_1}{L_1} \text{ and } r_2 = \frac{B_2}{L_2}.$$

Then, the volume of the asymmetric house from Figure 4 is equal to:

$$V = H[L_1 L_2 - L_1(1 - r_1)L_2(1 - r_2)] = H L_1 L_2 (r_1 + r_2 - r_1 r_2). \tag{14}$$

Therefore, height $H$ can be expressed as

$$H = \frac{V}{L_1 L_2 (r_1 + r_2 - r_1 r_2)}. \tag{15}$$

The external surface area of the asymmetric L-shaped house is given by

$$S = L_1 L_2 (r_1 + r_2 - r_1 r_2) + 2(L_1 + L_2)H = L_1 L_2 (r_1 + r_2 - r_1 r_2) + \frac{2V(L_1 + L_2)}{L_1 L_2 (r_1 + r_2 - r_1 r_2)}. \tag{16}$$

The first component in the formula above represent the area of the roof. The remaining part is equal to the area of the walls. Also in this case the area of the building in contact with the ground is not included in $S$.

The main goal of this section is to find the optimal L-shaped house geometry that minimizes the external surface area $S$. We will check different scenarios.

### 2.2.1 First Case: Fixed Volume $V$

In this scenario, the building volume $V$ is assumed fixed. However, similarly to the symmetric case, this formulation yields a minimum of the external surface area $S$ when

$L_1 = L_2 = B_1 = B_2$. In this situation, the building degenerates into a cuboid, and the characteristic L-shaped configuration is no longer preserved. Consequently, this solution is not acceptable within the context of L-shaped house design.

An analogous situation arises when the floor area $F$ and the building height $H$ are prescribed. In this case, the external envelope area is likewise minimized for $L_1 = L_2 = B_1 = B_2$, leading once again to a cuboid configuration. As a result, the defining L-shaped geometry is lost, rendering this solution unsuitable for the present design context.

### 2.2.2 Second Case: Fixed Wing Aspect Ratios $r_1$ and $r_2$

At the early stage of the design process, it is reasonable to treat the wing aspect ratios $r_1 = B_1/L_1$ and $r_2 = B_2/L_2$ as fixed parameters. In an asymmetric L-shaped house, the two wings typically serve different functional purposes, which imposes predefined requirements on room proportions, circulation, and minimum widths. These functional constraints are naturally expressed through prescribed aspect ratios rather than fully free geometric variables.

In addition, site conditions, zoning regulations, and setback requirements often restrict admissible wing proportions before detailed geometric optimization is undertaken. Structural and economic considerations, such as preferred spans, standardized components, and roof geometry, also favor specific length-to-width ratios. Fixing $r_1$ and $r_2$ therefore ensures that the optimization remains consistent with practical, regulatory, and constructional constraints, while preserving the intended asymmetric L-shaped character of the building.

In what follows we will find optimal shape (minimal envelope $S$) of the asymmetric L-shaped house, assuming that $r_1$ and $r_2$ as fixed parameters. As before, we assume that $V$ is also fixed.

Using Eq. (16), we get that the equations $\frac{\partial S}{\partial L_1} = 0$ and $\frac{\partial S}{\partial L_2} = 0$ take the form

$$L_2(r_1 + r_2 - r_1 r_2) + \frac{2V}{L_1 L_2 (r_1 + r_2 - r_1 r_2)} - \frac{2V(L_1 + L_2)}{L_1^2 L_2 (r_1 + r_2 - r_1 r_2)} = 0 \;,$$

$$L_1(r_1 + r_2 - r_1 r_2) + \frac{2V}{L_1 L_2 (r_1 + r_2 - r_1 r_2)} - \frac{2V(L_1 + L_2)}{L_1 L_2^2 (r_1 + r_2 - r_1 r_2)} = 0 \;.$$

After some tedious calculations, we get that the solution of the above system of equations, that defines the minimum of $S$, has the form

$$L_1^{\min} = L_2^{\min} = \left(\frac{2V}{(r_1 + r_2 - r_1 r_2)^2}\right)^{1/3}. \tag{17}$$

Consequently,

$$B_1^{\min} = r_1 L_1^{\min} = r_1 \left(\frac{2V}{(r_1 + r_2 - r_1 r_2)^2}\right)^{1/3}, \tag{18}$$

$$B_2^{\min} = r_2 L_2^{\min} = r_2 \left(\frac{2V}{(r_1 + r_2 - r_1 r_2)^2}\right)^{1/3}, \tag{19}$$

$$H^{\min} = \frac{V}{L_1^{\min} L_2^{\min}(r_1 + r_2 - r_1 r_2)} = \left(\frac{(r_1 + r_2 - r_1 r_2)V}{4}\right)^{1/3}, \tag{20}$$

$$S^{\min} = L_1^{\min} L_2^{\min}(r_1 + r_2 - r_1 r_2) + \frac{2V(L_1^{\min} + L_2^{\min})}{L_1^{\min} L_2^{\min}(r_1 + r_2 - r_1 r_2)} = 3\left(\frac{4V^2}{r_1 + r_2 - r_1 r_2}\right)^{1/3}. \tag{21}$$

The expressions given in (17)-(21) provide a practical framework for the design of an asymmetric L-shaped house that achieves the minimum external surface area $S_{\min}$ for a prescribed building volume $V$ and a fixed wing aspect ratios $r_1$ and $r_2$. Employing these relations allows designers to optimize the building geometry in a systematic manner, leading to a reduction in material usage and construction costs.

As a representative example, consider a building volume of $V = 300 m^3$ and aspect ratios $r_1 = 0{,}4, r_2 = 0{,}6$. Substituting these values into (17)-(21) yields the optimal parameters $L_1^{\min} = L_2^{\min} = 10{,}13 m, B_1^{\min} = 4{,}05 m$, $B_2^{\min} = 6{,}08 m, H^{\min} = 3{,}85 m$, and $S^{\min} = 233{,}86 m^2$.

From a practical standpoint, the obtained dimensions are well aligned with typical residential construction standards. They represent feasible design solutions that balance structural efficiency, economic viability, and enhanced energy performance. Figure 5 depicts the corresponding external surface area $S$ as a function of the wings lengths $L_1$ and $L_2$ The global minimum of $S$, corresponding to the optimal solution obtained analytically in Eq. (21), is indicated by a red marker. The figure therefore provides a graphical validation of the analysis.

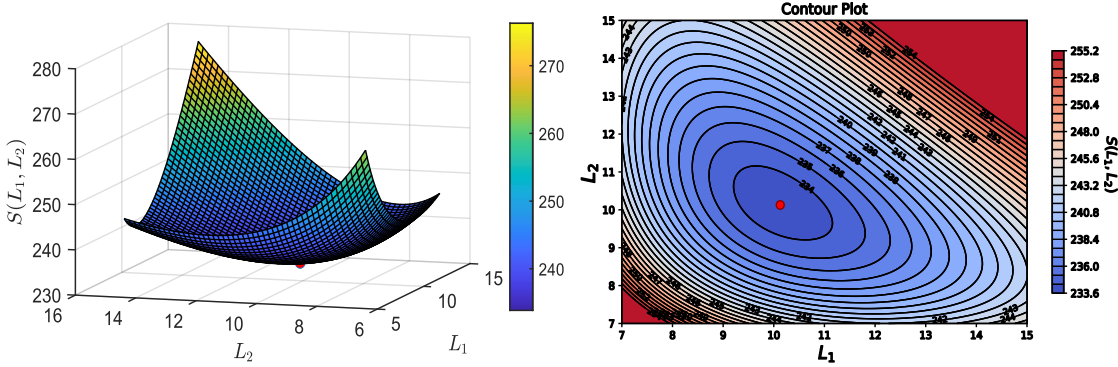

Fig. 5. Left panel: plot of the external surface $S(L_1, L_2)$. Right panel: the corresponding contour plot. The red dot indicates the location of the minimal surface obtained using formula (21). Parameters: $V = 300 m^3, r_1 = 0{,}4, r_2 = 0{,}6$.

*2.2.3 Third Case: Wing Aspect Ratios $r_1$ and $r_2$ from Fixed Ranges*

At the early stage of the design process, it is reasonable to assume that the wing aspect ratios $r_1$ and $r_2$ lie within a prescribed intervals $r_1 \in [a_1, b_1]$ and $r_2 \in [a_2, b_2]$, rather than being completely unrestricted. This assumption is primarily motivated by practical considerations arising from constraints that typically precede any detailed geometric optimization of the building envelope. The aspect ratios have a significant impact on the functional quality of each wing. In addition, site-specific factors, including plot geometry, setback regulations, easements, minimum boundary distances, and local planning rules, often further constrain the range of admissible building proportions.

For these reasons, we assume that the ratios are constrained to prescribed intervals: $r_1 \in [a_1, b_1]$ and $r_2 \in [a_2, b_2]$. As before, the volume $V$ is fixed. Under these assumptions, our objective is to determine the optimal house geometry that minimizes the external surface area $S$, as defined in (16).

To this end, we employ the method of Lagrange with KKT conditions, see the analysis in Sec. 2.1.3. and [26]. More precisely, we will find the optimal parameters $L_1^*, L_2^*, r_1^*$ and $r_2^*$ that minimize the external surface area $S$, in accordance with (16), The constrains are

$$a_1 \leq r_1 \leq b_1, a_2 \leq r_2 \leq b_2, L_1 > 0, L_2 > 0.$$

Here $0 < a_1 < b_1 < 1, 0 < a_2 < b_2 < 1$ and $V > 0$ are fixed constants. The corresponding Lagrangian is

$$\mathcal{L} = \frac{2V(L_1 + L_2)}{L_1 L_2 (r_1 + r_2 - r_1 r_2)} + L_1 L_2 (r_1 + r_2 - r_1 r_2) + \lambda_1 (a_1 - r_1) + \lambda_2 (r_1 - b_1)$$
$$+ \lambda_3 (a_2 - r_2) + \lambda_4 (r_2 - b_2).$$

To find the optimal values of $L_1, L_2, r_1$ and $r_2$ we apply the KKT conditions, cf. Sec. 2.1.3. Solving the system of KKT condition, after some tedious calculations we get that the optimal parameters minimizing $S$ are the following:

$$r_1^* = b_1, r_2^* = b_2, \tag{22}$$

$$L_1^* = L_2^* = \left( \frac{2V}{(b_1 + b_2 - b_1 b_2)^2} \right)^{1/3}, \tag{23}$$

and therefore

$$S^* = 3 \cdot \left( \frac{4V^2}{b_1 + b_2 - b_1 b_2} \right)^{1/3}. \tag{24}$$

The remaining optimal parameters are

$$B_1^* = b_1 L_1^* = b_1 \left( \frac{2V}{(b_1 + b_2 - b_1 b_2)^2} \right)^{1/3}, \tag{25}$$

$$B_2^* = b_2 L_2^* = b_2 \left( \frac{2V}{(b_1 + b_2 - b_1 b_2)^2} \right)^{1/3}, \tag{26}$$

$$H^* = \frac{V}{L_1^* L_2^* (b_1 + b_2 - b_1 b_2)} = \left( \frac{(b_1 + b_2 - b_1 b_2)V}{4} \right)^{1/3}. \tag{27}$$

The results summarized in (22)-(27) provide a practical foundation for the design of an L-shaped house that minimizes the external surface area for a prescribed building volume and a wing aspect ratios constrained to predefined intervals $r_1 \in [a_1, b_1]$ and $r_2 \in [a_2, b_2]$. In particular, they demonstrate that the optimal values of the aspect ratios within these intervals are attained at their upper bounds, that is, $r_1 = b_1$ and $r_2 = b_2$.

As an illustrative example, consider a building volume of $V = 200 \ m^3$ with the admissible aspect-ratio ranges defined by $a_1 = 0{,}3, b_1 = 0{,}5, a_2 = 0{,}2$ and $b_2 = 0{,}8$. Substituting these values into (22)-(27) yields the optimal parameters $r_1^* = 0{,}5, r_2^* = 0{,}8, L_1^* = L_2^* = 7{,}90 \ m, B_1^* = 3{,}95 \ m, B_2^* = 6{,}32 \ m, H^* = 3{,}95 \ m$, and $S^* = 168{,}69 \ m^2$.

From a practical standpoint, these dimensions are well aligned with typical residential construction standards. They represent feasible design solutions that balance structural efficiency, economic viability, and enhanced energy performance.

Figure 6 illustrates the plot of the minimal surface calculated for each $r_1 \in [a_1, b_1]$ and $r_2 \in [a_2, b_2]$ using formula (21). The global minimum of $S$ is indicated by a red marker.

This minimum is attained at the upper bounds of the admissible aspect-ratios intervals $r_1 = b_1$ and $r_2 = b_2$, confirming the analytical results (22)-(27).

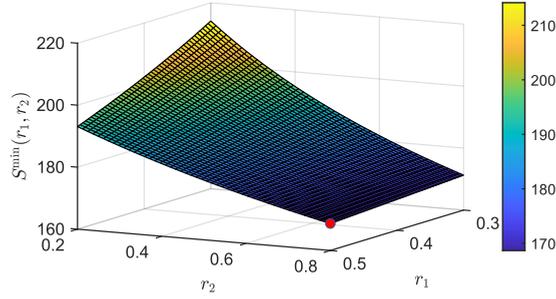

Figure 6: Plot of the minimal surface calculated for each $r_1 \in [a_1, b_1]$ and $r_2 \in [a_2, b_2]$ using formula (21). The red dot indicates the global minimum of $S$. Clearly, the global minimum is obtained for $r_1 = b_1$ and $r_2 = b_2$, in accordance with results (22)-(27). Parameters: $V = 200\ m^3, a_1 = 0{,}3, b_1 = 0{,}5, a_2 = 0{,}2$ and $b_2 = 0{,}8$.

### 2.2.4 Fourth Case: Fixed Height H and Wing Aspect Ratios $r_1$ and $r_2$

At the early design stage, it is reasonable to treat the building height as a predefined parameter, since it is typically constrained by factors that precede geometric optimization. Local zoning and planning regulations often impose strict limits on allowable building height, number of storeys, and roof geometry. In addition, functional requirements such as desired ceiling heights, storey-to-storey dimensions, and the inclusion of attic or technical spaces are usually specified by the investor at the outset.

So, now we assume that $H, r_1, r_2$ and $V$ are fixed. In such setting the building envelope area is equal to

$$S = \frac{V}{H} + 2H\left(L_1 + \frac{V}{L_1 H(r_1 + r_2 - r_1 r_2)}\right).$$

Minimum of this function is obtained for

$$L_1^{\min} = \sqrt{\frac{V}{H(r_1 + r_2 - r_1 r_2)}}. \tag{28}$$

The remaining optimal parameters are

$$L_2^{\min} = L_1^{\min}, B_1^{\min} = r_1 L_1^{\min}, B_2^{\min} = r_2 L_2^{\min}, S^{\min} = \frac{V}{H} + 4H\sqrt{\frac{V}{H(r_1 + r_2 - r_1 r_2)}}. \tag{29}$$

The expressions given in (28)-(29) establish a practical basis for the design of an asymmetric L-shaped house that attains the minimum external surface area $S_{\min}$ for prescribed $H, V, r_1$ and $r_2$. Making use of these relations enables a systematic optimization of the building geometry, which in turn contributes to reduced material consumption and lower construction costs.

## 3. CASE STUDIES

To emphasize the practical significance of the theoretical optimality results derived in the previous section, we consider an application in the form of two case studies. In particular, we analyze two L-shaped house designs currently offered by architectural firms, denoted by A and B (see Fig. 7). These buildings were intentionally selected to represent a diverse range of geometric shapes, proportions, and sizes. The objective is to illustrate how the derived formulas can be applied to the analysis of real-world structures with varying dimensions and configurations. For each house, both the visual characteristics and the key geometric parameters were extracted using data available on project websites. This methodology ensures that the case studies are not only illustrative but also reproducible, as it can be readily extended to other architectural examples.

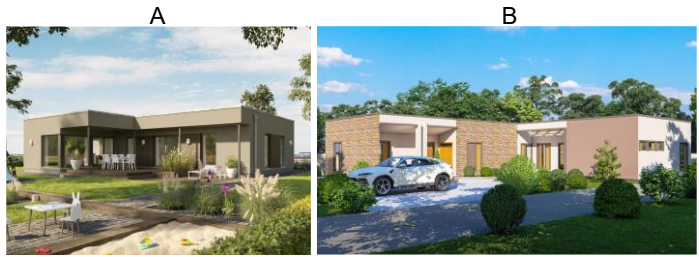

Fig. 7. Images of two L-shaped houses analyzed as case studies using the results obtained in previous sections. Design and image sources: www.hausbaudirekt.de/haus/solution-110-v4-2-2 (panel A); www.maramani.com/products/l-shaped-4-bedroom-house-1d-14313 (panel B).

**Analysis of House A.** The parameters of House A are as follows: $L_1 = 13,7m, L_2 = 14,9m, B_1 = 8,7m, B_2 = 4,6m, H = 3,6m$. This gives $r_1 = 0,63, r_2 = 0,31$ V $= 549,5\ m^3$ and $S = 358,5\ m^2$.

Let us apply the optimality results from (17)-(21). We get the following optimal parameters of house A : $L_1^{\min} = L_2^{\min} = 12,5m, B_1^{\min} = 7,9m, B_2^{\min} = 3,9m, H^{\min} = 4,6m$ and $S^{min} = 351{,}9 m^2$ Th above optimization procedure would reduce the external surface by $6{,}6\ m^2$.

Application of the optimality results from (28)-(29) gives the external surface almost identical to the original one. From this point of view, the original design is very close to optimal. This result is expected due to very limited degrees of freedom in (28)-(29).

**Analysis of House B.** The parameters of House B are the following: $L_1 = 22; m, L_2 = 19{,}5m, B_1 = 8{,}8m,\ B_2 = 8{,}5m, H = 4{,}3m$. This gives $r_1 = 0{,}4, r_2 = 0{,}43 V = 1220{,}3 m^3$ and $S = 640{,}7 m^2$.

We apply the optimality results from (17)-(21). We obtain the following optimal parameters of house B: $L_1^{min} = L_2^{min} = 17{,}7m, B_1^{min} = 7{,}1m, B_2^{min} = 7{,}7m, H^{min} = 5{,}8m$ and $S^{min} = 624{,}1 m^2$ Th above optimization procedure would reduce the external surface by $16{,}6\ m^2$. It would result in reasonable savings in construction costs and energy consumption of the house.

Note that the application of the optimality results from (28)-(29) gives the external surface smaller than the original one by about $1 m^2$. From this point of view, the original design is close to optimal. This result is expected due to very limited degrees of freedom in (28)-(29).

## 4. CONCLUSIONS

This paper has presented a systematic analytical study of the geometric optimization of L-shaped residential buildings with the objective of minimizing the external surface area for a prescribed building volume. Since the external surface area is directly related to material consumption, heat losses, and operational energy demand, the obtained results provide a clear geometric basis for improving both energy efficiency and sustainability at the early design stage.

For symmetric L-shaped houses, several optimization scenarios were investigated. It was shown that when only the volume or the floor area is fixed, the minimum external surface area is attained for a degenerate cuboid configuration, leading to the loss of the L-shaped form. This observation highlights the necessity of introducing additional geometric constraints in order to preserve the intended building topology. By fixing the wing aspect ratio, explicit formulas for the optimal wing dimensions, height, and minimal envelope area were derived. These closed-form expressions allow designers to directly determine optimal proportions that reduce envelope area while maintaining functional and architectural requirements. When the wing aspect ratio was constrained to a prescribed interval, the constrained optimization analysis based on KKT conditions demonstrated that the optimal solution is attained at the boundary of the admissible range, providing a clear design guideline for selecting feasible proportions. An analogous optimization framework was developed for asymmetric L-shaped houses, which offer greater flexibility and are particularly relevant for irregular plots and context-driven designs. For fixed wing aspect ratios, the optimal configuration was shown to correspond to equal wing lengths, despite differing widths, and explicit expressions for the optimal geometry and minimal external surface area were obtained. When the aspect ratios were constrained to admissible intervals, the analysis revealed that the minimum envelope area is achieved at the upper bounds of these intervals. Additional scenarios with fixed building height further illustrated how regulatory and functional constraints can be incorporated into the optimization process without sacrificing analytical tractability.

The practical relevance of the theoretical results was demonstrated through real-world case studies of existing L-shaped house designs. In both cases, the application of the derived formulas led to a reduction of the external surface area or confirmed that the original designs were already close to optimal. These examples illustrate that even modest geometric adjustments, guided by analytical optimization, can lead to tangible savings in construction materials and energy demand.

Overall, the presented results show that geometric optimization of L-shaped buildings can be effectively addressed using relatively simple analytical tools, yielding explicit design formulas that are directly applicable in architectural practice. The proposed approach supports informed early-stage decision-making, enabling designers to balance functional requirements, regulatory constraints, architectural intent, and energy efficiency. Future research may extend this framework to include additional factors such as thermal zoning, façade orientation, window-to-wall ratios, or multi-storey configurations, further strengthening the link between geometric form and sustainable building performance.